\documentclass[microtype]{gtpart}     

\usepackage{dsfont}
\usepackage{hyperref}
\usepackage{mathrsfs}
\usepackage{graphicx}
\usepackage{mathdesign}
\usepackage{float}
\usepackage{tikz}
\usetikzlibrary{matrix,arrows}
\usepackage{epstopdf}

\title[The topology of nilpotent representations in reductive groups]{The topology of nilpotent representations in reductive groups and their maximal compact subgroups}

\author{Maxime Bergeron}
\givenname{Maxime}
\surname{Bergeron}
\address{Department of Mathematics, The University of British Columbia, Room 121 - 1984 Mathematics Road, V6T 1Z2, Vancouver, Canada}
\email{mbergeron@math.ubc.ca}
\urladdr{http://www.math.ubc.ca/~mbergeron/}


\arxivreference{	arXiv:1310.5109 [math.GR]}
\arxivpassword{2ruid}

\volumenumber{}
\issuenumber{}
\publicationyear{}
\papernumber{}
\startpage{}
\endpage{}
\doi{}
\MR{}
\Zbl{}
\received{}
\revised{}
\accepted{}
\published{}
\publishedonline{}
\proposed{}
\seconded{}
\corresponding{}
\editor{}
\version{}

\newtheorem{theorem}{Theorem}[section]
\newtheorem{lemma}[theorem]{Lemma}
\newtheorem{proposition}[theorem]{Proposition}
\newtheorem{corollary}[theorem]{Corollary}

\newtheorem*{namedtheorem}{\theoremname}
\newcommand{\theoremname}{testing}
\newenvironment{named}[1]{\renewcommand{\theoremname}{#1}\begin{namedtheorem}}{\end{namedtheorem}}

\theoremstyle{remark}

\newenvironment{definition}[1][Definition]{\begin{trivlist}
\item[\hskip \labelsep {\bfseries #1}]}{\end{trivlist}}

\newenvironment{example}[1][Example]{\begin{trivlist} 
\item[\hskip \labelsep {\bfseries #1}]}{\end{trivlist}}

\newenvironment{remark}[1][Remark]{\begin{trivlist}
\item[\hskip \labelsep {\bfseries #1}]}{\end{trivlist}}

\newcommand{\bb}[1]{\mathbb{#1}}
\newcommand{\fk}[1]{\mathfrak{#1}}

\newcommand{\raw}{\rightarrow}

\newcommand{\inv}{^{-1}}

\DeclareMathOperator{\Hom}{Hom}	
\DeclareMathOperator{\sln}{SL_n\bb{C}}

\DeclareMathOperator{\GL}{GL}
\DeclareMathOperator{\SL}{SL}
\DeclareMathOperator{\SU}{SU}

\DeclareMathOperator{\On}{O_n}
\DeclareMathOperator{\One}{O^+_{n,1}}
\DeclareMathOperator{\U}{U}
\DeclareMathOperator{\spn}{Sp_n\bb{C}}
\DeclareMathOperator{\Aut}{Aut}
\DeclareMathOperator{\Id}{Id}

\begin{document}

\begin{abstract}   
Let $G$ be a complex reductive linear algebraic group and let $K\subset G$ be a maximal compact subgroup. Given a nilpotent group $\Gamma$ generated by $r$ elements, we consider the representation spaces $\Hom(\kern 0.05em\Gamma,G\kern 0.05em)$ and $\Hom(\kern 0.05em\Gamma,K\kern 0.05em)$ with the natural topology induced from an embedding into $G^r$ and $K^r$ respectively. The goal of this paper is to prove that there is a strong deformation retraction of $\Hom(\kern 0.05em\Gamma,G\kern 0.05em)$ onto  $\Hom(\kern 0.05em\Gamma,K\kern 0.05em)$. We also obtain a strong deformation retraction of the geometric invariant theory quotient $\Hom(\kern 0.05em\Gamma,G\kern 0.05em)/\!\!/G$ onto the ordinary quotient $\Hom(\kern 0.05em\Gamma,K\kern 0.05em)/K$.
\end{abstract}

\maketitle


\section{Introduction}

Let $\Gamma$ be a group generated by $r$ elements and let $G$ be a complex linear algebraic group. Since a representation of $\Gamma$ in $G$ is uniquely determined by the image of a generating set, the space $\Hom(\kern 0.05em\Gamma,G\kern 0.05em)$ of homomorphisms $\rho\negthinspace:\Gamma\raw G$ can be realized as an affine algebraic set carved out of $G^r$ by the relations of $\Gamma$. It is well known and otherwise easy to see that this geometric structure is independent of the chosen presentation of $\Gamma$.
As a complex variety, $\Hom(\kern 0.05em\Gamma,G\kern 0.05em)$ admits a natural Hausdorff topology obtained from an embedding into affine space. If $K$ is a maximal compact subgroup of $G$, we can then endow $\Hom(\kern 0.05em\Gamma,K\kern 0.05em)\subset\Hom(\kern 0.05em\Gamma,G\kern 0.05em)$ with the subspace Hausdorff topology. Although in general these topological spaces may be quite different, in this paper we show that when $\Gamma$ is nilpotent and $G$ is reductive they are actually homotopy equivalent: 

\begin{named}{Theorem I}\hypertarget{Theorem I}{}
Let $\Gamma$ be a finitely generated nilpotent group and let $G$ be the group of complex or real points of a (possibly disconnected) reductive linear algebraic group, defined over $\bb{R}$ in the latter case. If $\,K$ is any maximal compact subgroup of $\,G$, then there is a ($K$--\,equivariant) strong deformation retraction of $\,\Hom(\kern 0.05em\Gamma,G\kern 0.05em)$ onto $\Hom(\kern 0.05em\Gamma,K\kern 0.05em)$. In particular, $\Hom(\kern 0.05em\Gamma,G\kern 0.05em)$ and $\,\Hom(\kern 0.05em\Gamma,K\kern 0.05em)$ are homotopy equivalent.
\end{named}

The assumption that $G$ be reductive in \hyperlink{Theorem I}{Theorem I} is necessary. Combining ideas of Mal'cev and Dyer, one can produce a (nilpotent) lattice $\Gamma$ in a unipotent Lie group $U$ for which the space $\Hom(\kern 0.05em\Gamma,U\kern 0.05em)$ is disconnected. Viewing $U$ as an algebraic group and keeping in mind that a maximal compact subgroup $K\subset U$ is necessarily trivial, we see that the inclusion $\Hom(\kern 0.05em\Gamma,K\kern 0.05em)\hookrightarrow \Hom(\kern 0.05em\Gamma,U\kern 0.05em)$ is not even a bijection of $\pi_0$. A complete discussion of this example will be given in Appendix \ref{appendix}. Similarly, \hyperlink{Theorem I}{Theorem I} does not hold for arbitrary $\Gamma$: if $\Gamma\subset \SL_2\!\bb{C}$  is a cocompact lattice, then Mostow Rigidity shows that the tautological representation $\Gamma\hookrightarrow \SL_2\!\bb{C}$ lies in a component of $\Hom(\kern 0.05em\Gamma,\SL_2\!\bb{C}\kern 0.05em)$ which is disjoint from $\Hom(\kern 0.05em\Gamma,\SU_2\kern 0.05em)$. Once again the inclusion $\Hom(\kern 0.05em\Gamma,K\kern 0.05em)\hookrightarrow\Hom(\kern 0.05em\Gamma,G\kern 0.05em)$ isn't  a bijection of $\pi_0$.

There is also an interesting contrast between the real and complex cases of \hyperlink{Theorem I}{Theorem I}.
For instance, while nilpotent subgroups of $\One\!\bb{R}$ (the group of isometries of hyperbolic $n$-space) are all virtually abelian, this is far from being true in $\GL_n\!\bb{R}$. 
Nevertheless, $\On\!\bb{R}$ is a maximal compact subgroup of both $\GL_n\!\bb{R}$ and $\One\!\bb{R}$ so \hyperlink{Theorem I}{Theorem I} shows that $\Hom(\kern 0.05 em\Gamma,\One\!\bb{R})$ and $\Hom(\kern 0.05 em\Gamma,\GL_n\!\bb{R})$ are homotopy equivalent. This kind of example does not occur in the complex case because complex reductive linear algebraic groups are determined by their maximal compact subgroups.

It should be mentioned at this point that \hyperlink{Theorem I}{Theorem I} has two predecessors. When $\Gamma$ is an abelian group, it is due to Pettet and Souto \cite{pettet2011commuting} and when $\Gamma$ is an expanding nilpotent group (c.f. Appendix \ref{appendix}) it is is due to Lior Silberman and Juan Souto. While these authors used homotopy-theoretic methods, in this paper we rely on algebraic geometry instead. In particular, even in the case where $\Gamma=\bb{Z}^r$ our proof is completely new. One of the advantages of this approach, coming from geometric invariant theory, is that the corresponding result for the character variety $\Hom(\kern 0.05em\Gamma,G\kern 0.05em)/\!\!/G$ follows from the proof of \hyperlink{Theorem I}{Theorem I}:

\begin{named}{Theorem II}\hypertarget{Theorem II}{}
Let $\,\Gamma$ be a finitely generated nilpotent group and let $G$ be the group of complex or real points of a (possibly disconnected) reductive linear algebraic group, defined over $\bb{R}$ in the latter case. If $K\subset G$ is any maximal compact subgroup, then there is a strong deformation retraction of the character variety $\Hom(\kern 0.05em\Gamma,G\kern 0.05em)/\!\!/G$ onto $\Hom(\kern 0.05em\Gamma,K\kern 0.05em)/K$.
\end{named}
When $\Gamma$ is abelian, \hyperlink{Theorem II}{Theorem II} is a recent result of Florentino and Lawton \cite{florentino2014topology} (following up on their previous work for free groups in \cite{florentino2009topology}) and we refer the reader to their paper for various applications.
The main tools from geometric invariant theory used to prove \hyperlink{Theorem I}{Theorem I} and \hyperlink{Theorem II}{Theorem II} are the Kempf--Ness Theorem and the Neeman--Schwarz Theorem. More precisely, for $\Gamma$ and $G$ as in \hyperlink{Theorem I}{Theorem I}, we embed $\Hom(\kern 0.05em\Gamma,G\kern 0.05em)$ into the complex vector space $\oplus_{j=1}^r M_n\bb{C}$ equipped with the Frobenius norm. Following Richardson, we then denote by  $\mathcal{M}$ the set of representations that are of minimal norm within their $G$--orbit, the so-called \emph{Kempf--Ness set} of $\Hom(\kern 0.05em\Gamma,G\kern 0.05em)$. Our main technical result characterizes this set when $\Gamma$ is nilpotent:
\begin{named}{Theorem III}\hypertarget{Theorem III}{}
Let $G\subset \sln$ be a complex reductive linear algebraic group and suppose that $K=G\cap \SU_n\subset G$ is a maximal compact subgroup. If $\,\Gamma$ is a finitely generated nilpotent group and $\mathcal{M}$ denotes the Kempf--Ness set of  the $G$--subvariety $\Hom(\kern 0.05em\Gamma,G\kern 0.05em)\subset \oplus_{j=1}^r M_n\bb{C}$ equipped with the Frobenius norm, then $$
\mathcal{M}=\{\rho\in\Hom(\kern 0.05em\Gamma,G\kern 0.05em):\rho(\Gamma)\text{ consists of normal matrices}\}.
$$
\end{named} 
Here, it is worth noting that the nilpotency assumption on $\Gamma$ is necessary; indeed, \hyperlink{Theorem III}{Theorem III} is false for the solvable group $\bb{Z}\rtimes \bb{Z}/4\bb{Z}$ (c.f. Section \ref{the Kempf--Ness set}). Since the Neeman--Schwarz Theorem provides us with a retraction of $\Hom(\kern 0.05em\Gamma,G\kern 0.05em)$ onto $\mathcal{M}$, in order to deduce \hyperlink{Theorem I}{Theorem I} from \hyperlink{Theorem III}{Theorem III} it remains to deform representations of normal matrices into representations of unitary matrices. This is accomplished by scaling eigenvalues. Finally, \hyperlink{Theorem II}{Theorem II} is a consequence of \hyperlink{Theorem III}{Theorem III} after applying a corollary of the Kempf--Ness Theorem.

\subsection{Historical remarks and applications}
In coarsest terms, this paper is concerned with the geometric classification of representations of finitely generated groups. Classically, this subject finds its roots in Poincar\'e's work  on monodromy groups of linear homogeneous equations and geometric invariant theory. 
More recently, the study of representation varieties has had impacts in a variety of contexts and we mention here but a few.

From a first point of view, the geometry of representation varieties can be used to deduce algebraic information about the representation theory of a group. For instance, 
 Lubotzky and Magid \cite{lubotzky1985varieties} obtained such information for the $\GL_n\!\bb{C}$ representations of a finitely generated group $\Gamma$ by using Weil's results on the Zariski tangent space of $\Hom(\kern 0.05 em\Gamma,\GL_n\!\bb{C}\kern 0.05 em)$. Incidentally, these techniques were most effective when $\Gamma$ was assumed to be nilpotent, partly because this ensured a special vanishing property of its first group cohomology.

 One can also study representation varieties from a purely differential geometric point of view. Here, when $\Gamma$ is the fundamental group of a smooth manifold $M$, $\Hom(\kern 0.05em\Gamma,K\kern 0.05em)$ can be identified with the space of  pointed flat connections on principal $K$--bundles over $M$.  These spaces lie at the intersection of various fields including gauge theory and symplectic geometry as illustrated for surfaces in Jeffrey's survey \cite{jeffrey2005flat}. More precisely,
there is a natural conjugation action of $K$ on $\Hom(\kern 0.05em\Gamma,K\kern 0.05em)$ and the quotient $\Hom(\kern 0.05em\Gamma,K\kern 0.05em)/K$ corresponds to the moduli space of flat-connections on principal $K$--bundles over $M$. In the case where $M$ is compact and K{\"a}hler, the geometric invariant theory quotient $\Hom(\kern 0.05em\Gamma,G\kern 0.05em)/\!\!/G$  corresponds to the moduli space of polystable $G$--bundles over $M$ via the Narasimhan--Seshadri Theorem (see Narasimhan--Seshadri \cite{narasimhan1964holomorphic} and  Simpson \cite{simpson1994moduli}). By and large, these so-called character varieties have attracted a lot of attention lately as shown in Sikora's survey  
\cite[Section 11]{sikora2012character}
 and the references therein.

More topologically, the work in supersymmetric Yang--Mills theory and mirror symmetry of Witten in \cite{witten1982constraints} and \cite{witten1998toroidal} sparked an interest in the connected components of free abelian representation varieties as seen in Kac--Smilga \cite{kac2000vacuum} and Borel--Friedman--Morgan \cite{borel2002almost}. Here, there has been much recent development in the study of higher topological invariants stemming from work of \'Adem and Cohen \cite{adem2007commuting} in the compact case. Prior to Pettet and Souto's result \cite{pettet2011commuting}, most known topological invariants concerned spaces of representations into compact groups; it was 
 their deformation retraction which allowed many of these results to be extended to representations into reductive groups.  We briefly indicate how this strategy works in the nilpotent case,  referring the reader to \cite{pettet2011commuting} for more applications of this type. 
 
Suppose for the sake of concreteness that $\Gamma$ is a torsion free nilpotent group and $G$ is a complex reductive linear algebraic group with a given maximal compact subgroup $K\subset G$. Denote by $\Hom(\kern 0.05em\Gamma,K\kern 0.05em)_\mathds{1}$ the connected component of $\Hom(\kern 0.05em\Gamma,K\kern 0.05em)$ containing the trivial representations.
In \cite{gomez2012fundamental}, G\'omez, Pettet and Souto   compute that $\pi_1(\,\Hom(\kern 0.05 em\bb{Z}^r,K\kern 0.05 em)_\mathds{1})\cong \pi_1(K)^r$ and use the retraction constructed in \cite{pettet2011commuting} to conclude that $\pi_1(\,\Hom(\kern 0.05 em\bb{Z}^r,G\kern 0.05 em)_\mathds{1})\cong \pi_1(G)^r$. 
 It is not hard to see in our case that  $\Hom(\kern 0.05em\Gamma,K\kern 0.05em)_\mathds{1}$ coincides with the representations factoring through the abelianization of $\Gamma$ lying in the component of the trivial representation. These observations can be combined with \hyperlink{Theorem I}{Theorem I} to compute that $\pi_1(\,\Hom(\kern 0.05em\Gamma,G\kern 0.05em)_\mathds{1})\cong \pi_1(G)^{\,\mathrm{rank}\, H_1(\Gamma;\bb{Z})}.$ The topology of $\Hom(\kern 0.05em\Gamma,K\kern 0.05em)$ will be further analyzed in subsequent papers.

\subsection{Outline of the paper}

We begin in Section \ref{algebgroups} by establishing some notation and refreshing the reader's memory with some basic facts about algebraic groups. Then, in Section \ref{algebraic actions}, we introduce the main technical tools from Kempf--Ness theory that we shall need. After these preliminaries, we begin working with representation varieties in Section \ref{the representation variety} where we prove \hyperlink{Theorem III}{Theorem III}. Then, in Section \ref{the scaling operation}, we show  how eigenvalues can be scaled to complete the proof of \hyperlink{Theorem I}{Theorem I} in the complex case. Finally, in Section \ref{real}, we show how the proof carries over to the real case and we prove \hyperlink{Theorem II}{Theorem II}.

\subsection{Acknowledgements} 

The author wishes to thank Lior Silberman for crucial conversations about algebraic groups and nilpotent groups. He would also like to thank Alejandro \'Adem, Galo Higuera Rojo, Justin Martel, Alexandra Pettet, Tali Pinsky and Zinovy Reichstein for carefully listening to his various rambles.
Most importantly, the author wishes to thank Juan Souto for introducing him to the topology of representation varieties and suggesting that he should explore the algebro-geometric approach. This paper would not have been possible without the many insightful conversations they had together.

\section{Algebraic groups}
In this preliminary section, we refresh the reader's memory with some basic facts about linear algebraic groups, referring to Borel's book \cite{borel1991linear} for the details. In his own words: 
\emph{According to one's taste about naturality and algebraic geometry, it is possible to give several definitions of linear algebraic groups.} Here, we will favour a concrete and slightly pedestrian approach since it is all that we shall need. There is only one possibly nonstandard result, Lemma \ref{componentnilp}.

\label{algebgroups}\subsection{A crash course}
An affine algebraic group is an affine  variety endowed with a group law for which the group operations are morphisms of varieties, i.e., polynomial maps. For our purposes, an algebraic group $G$ shall mean the group of complex (or real) points of an affine algebraic group (defined over $\bb{R}$ in the latter case) and an affine variety shall mean an affine algebraic set, i.e., the zero locus of a family of complex polynomials. It turns out that affine algebraic groups are \emph{linear} so we will always identify them with a Zariski closed subgroup of $\sln$.

\begin{remark}
We shall have the occasion to consider two different topologies on varieties and their subsets, the classical Hausdorff topology and the Zariski topology. 
 Unless we specify otherwise, all references to topological concepts will refer to the Hausdorff topology. For instance, a \emph{closed set} is closed in the Hausdorff topology and a \emph{Zariski closed} set is closed in the Zariski topology.
 \end{remark}

Let $G\subset \sln$ be an algebraic group.
\begin{enumerate}
\item We say that a subgroup $L$ is an \emph{algebraic subgroup} of $G$ if it is   Zariski closed.  Since an arbitrary subgroup $H$ of $G$ isn't necessarily algebraic, we often pass from $H$ to its algebraic closure $\mathscr{A}(H)$.  It turns out that $\mathscr{A}(H)$ is simultaneously the Zariski closure of $H$ in $G$ and the intersection of all algebraic subgroups of $G$ containing $H$. This same procedure can be applied when $H$ is a subset of $G$.

\item We say that $G$ is \emph{connected} if it is connected in the Zariski topology. Much of the behaviour of $G$ is governed by the connected component of its identity element which we denote by $G^o$. This is always a normal  algebraic subgroup of finite index in $G$. The following Theorem due to Mostow \cite{mostow1955self} 
 indicates one of the many r\^oles played by $G^o$; it shall be used several times in this paper.
\begin{theorem}[Conjugacy Theorem]\label{conjugacy}
An algebraic group $G$ contains a maximal compact subgroup and all such subgroups are conjugate by elements of $G^o$.
 \end{theorem}
\item Lastly,  since $G$ is a smooth variety, it has a well defined tangent space at the identity that we denote by $\fk{g}$. The \emph{Lie algebra} of $G$ is the set $\fk{g}$ endowed with the usual Lie algebra structure and bracket operation on derivations. 

\end{enumerate}

\subsection{Jordan decomposition}
Let $V$ be a finite dimensional complex vector space and recall the following elementary definitions from linear algebra; an endomorphism $\sigma$ of $V$ is said to be:
\begin{enumerate}
\item \emph{Nilpotent} if $\sigma^n=0$ for some $n\in\bb{N}$.
\item \emph{Unipotent} if $\sigma-\Id$ is nilpotent.
\item \emph{Semisimple} if $V$ is spanned by eigenvectors of $\sigma$.
\end{enumerate}
If $G\subset\sln$ is an algebraic group, we say that $g\in G$ is semisimple (resp. unipotent) if it is semisimple (resp. unipotent) as an endomorphism of $\bb{C}^n$. One can show that this definition does not depend on the chosen embedding of $G$ in $\sln$. In fact, we have the following useful theorem:

\begin{named}{Jordan Decomposition Theorem}
Let $G$ be a linear algebraic group.
\begin{enumerate}
\item If $g\in G$ then there are unique elements $g_s$ and $g_u$ in $G$ so that $g_s$ is semisimple, $g_u$ is unipotent and $g=g_sg_u=g_ug_s$. 
\item If $\varphi:G\raw G'$ is a morphism of algebraic groups, then $\varphi(g)_s=\varphi(g_s)$ and $\varphi(g)_u=\varphi(g_U\kern 0.05 em)$.
\end{enumerate}
\end{named}
We denote the set of all unipotent (resp. semisimple) elements of $G$ by $G_u$ (resp. $G_s$). Although the set $G_s$ is seldom Zariski closed, $G_u$ is always a Zariski closed subset of $G$ contained in $G^o$. In fact, we say that an algebraic group $G$ is \emph{unipotent} if $G=G_u$. We can now define the class of algebraic groups that we will be most concerned with in this paper:
\begin{definition}
The \emph{unipotent radical} of an algebraic group $G$ is its maximal connected normal unipotent subgroup. An algebraic group $G$ is said to be \emph{reductive} if its unipotent radical is trivial.
\end{definition}
\begin{example}
The classical groups are all reductive, e.g., $\sln$, $\GL_n\!\bb{C}$ and $\spn$.
\end{example}

\subsection{Nilpotent algebraic groups}\label{nilpred}
Let $A$ and $B$ be any subgroups of $G$ and denote by $[A,B]$ the commutator subgroup of $G$ generated by elements of the form $aba\inv b\inv$, $a\in A$ and $b\in B$. The \emph{lower central series} of $G$ is defined inductively by the rule
$$G=:G_{(0)}\unrhd[G,G_{(0)}]=:G_{(1)}\unrhd[G,G_{(1)}]=:G_{(2)}\unrhd\ldots\unrhd[G,G_{(n)}]=:G_{(n+1)}\unrhd\ldots$$  and one says that $G$ is  \emph{nilpotent} if for some $n\geq 0$ we have  $G_{(n)}=\{e\}$. This definition makes sense in the context of algebraic groups since all of the groups in the lower central series are algebraic, i.e., Zariski closed. In fact, we have the following classification result which may be found in Borel \cite[III.10.6]{borel1991linear}: 
\begin{proposition}\label{nilpotent decomposition}
Let $G$ be an algebraic group and recall that we denote its subset of semisimple (resp. unipotent) elements by $G_s$ (resp. $G_u$). If $G$ is connected and nilpotent, then $G_s$ and $G_u$ are algebraic subgroups of $G$ and $G\cong G_s\times G_u$. 
\end{proposition}

We end this subsection with a structural result that will be useful later on. Although it may be known to the experts, we provide a proof since we were unable to find it in the literature.

\begin{lemma}\label{componentnilp}
If $N$ is a (possibly disconnected) nilpotent reductive algebraic group, then 
\begin{enumerate}

\item $N$ consists of semisimple elements.
\item $N$ has a unique maximal compact subgroup.
\item The connected component $N^o$ is contained in the centre of $N$.

\end{enumerate}
\end{lemma}
The following argument was shown to us by Lior Silberman.
\begin{proof}
Recall from Proposition \ref{nilpotent decomposition} that the connected component of a nilpotent algebraic group admits a direct product decomposition $N^o=N^o_s\times N^o_u$. Since $N^o$ is reductive, $N^o_u$ must be trivial and $N^o=N^o_s$ consists of semisimple elements. The proof of claim $(1)$ is then completed by observing that $N_u\subset N^o$. Moreover, assuming for the moment that $N^o$ is central in $N$, claim $(2)$ follows at once from the Conjugacy Theorem \ref{conjugacy}.

 It remains to prove that $N^o$ is central in $N$. To do so, it suffices to show that the adjoint action \cite[Chapter I.3.13]{borel1991linear} of $N$ on its Lie algebra $\mathfrak{n}$ is trivial. Seeking a contradiction, suppose that $Ad(g)$ acts nontrivially on $\mathfrak{n}$ for some $g\in N$. Since $g$ is semisimple, the Jordan Decomposition Theorem ensures that $Ad(g)$ is semisimple and must therefore have an eigenvalue $\lambda\neq 1$. If $X\in\mathfrak{n}$ is an eigenvector of $Ad(g)$ corresponding to $\lambda$, we have that
  $$  [g,exp(tX)]=ge^{tX}g\inv e^{-tX}=e^{Ad(g)(tX)-tX}=e^{t(\lambda-1)X}.$$ As such, the iterated commutators $[g,[g,[g,\ldots,[g,[g,e^X]]\ldots ]]]$ are all non-trivial and the lower central series of $N$ does not terminate in finitely many steps. This is a contradiction.
\end{proof}

\subsection{Complex reductive algebraic groups}
\label{complexreductive}
One of the advantages of working over $\bb{C}$ is the strong relation between reductive groups and their maximal compact subgroups. We illustrate this by following the treatment in Schwarz \cite{schwarz1989topology}.

Given a compact Lie group $K$, the Peter-Weyl Theorem provides a faithful embedding $K\hookrightarrow \GL_n\!\bb{R}$ for some $n$. Identifying $K$ with its image realizes it as a real algebraic subgroup of $\GL_n\!\bb{R}$. We then define the \emph{complexification} $G:=K_\bb{C}$ to be the vanishing locus in $\GL_n\!\bb{C}$ of the ideal defining $K$. The group $G$ is a complex algebraic group which is independent  up to isomorphism of the embedding provided by the Peter-Weyl Theorem. The following properties of complexifications of compact Lie groups are well known (see, for instance, Onishchik and Vinberg \cite[Chapter 5]{onishchik1990lie}):
\begin{enumerate}
\item The Lie algebra $\fk{g}$ of $G$ has a natural splitting in terms of the Lie algebra $\fk{k}$ of $K$: $$\fk{g}=\fk{k}\oplus i\fk{k}.$$
\item The group $K$ is a maximal compact subgroup of $G$. Moreover, it is Zariski dense in $G$, i.e., $\mathscr{A}(K)=G$.
\item $G=K\cdot exp(i\fk{k})\cong K\times exp(i\fk{k}) $ where $\cong$ is a diffeomorphism. (Contracting $exp(i\fk{k})$ is one way to see that $K\hookrightarrow G$ is a homotopy equivalence.)
\end{enumerate}
\begin{example}
Viewing $\SU_n$ as a compact Lie group, we can realize it as a real algebraic subgroup of $\GL_{2n}\!\bb{R}$. In this case, the complexification of $\SU_n$ is isomorphic to $\SL_n\!\bb{C}\subset \GL_{2n}\!\bb{C}$. \end{example}
\begin{remark}
The decomposition $G=K\!\cdot\! exp(i\fk{k})$ is often called the \emph{polar decomposition}. Indeed, for any linear realization $G\subset\sln$, it is a Theorem of Mostow \cite{mostow1955self} that for some $m\in\sln$ the corresponding polar decomposition of $mGm\inv\subset\sln$ coincides with the usual polar decomposition in $\sln$. We shall use this decomposition several times in the proof of our main theorem.
\end{remark}

As one might suspect, the following theorem holds \cite[Chapter 5, Section 2]{onishchik1990lie}:
\begin{theorem}
A complex linear algebraic group is reductive if and only if it is the complexification of a compact Lie group.
\end{theorem}

\section{Algebraic actions}\label{algebraic actions}
Let $V$ be a complex vector space and suppose that the complex reductive algebraic group $G$ is contained in $\GL(V)$.  The goal of this preliminary section is to study the algebraic and topological structures of a variety $X\subset V$ that is stable under the action of $G\subset\GL(V)$.  This will lead us to the main known result used in our proof: the Neeman--Schwarz Theorem.

\label{actions}
In keeping with the general spirit of this paper, recall that we favour a concrete approach to algebraic groups. 
 As such, we also adopt the following concrete point of view for algebraic actions.
A \emph{G--variety} is a variety $X\subset V$ acted upon algebraically by $G\subset \GL(V)$.
 Recall that the \emph{orbit}  of $x\in X$, denoted $G\cdot x$, is the set of all $g\cdot x$ with $g\in G$.
 A \emph{G--subvariety} of $X$ is then a $G$--variety $Y\subset X\subset V$ such that $G\cdot y\subset Y$ for every $y\in Y$. If $X$ and $Y$ are $G$--varieties, a morphism $\alpha:X\raw Y$ is $G$\emph{--\,equivariant} if $\alpha(g\cdot x)=g\cdot \alpha (x)$.

 \begin{remark}
 There is no real loss of generality in our concrete approach to $G$--varieties. Indeed, any affine variety endowed with a reductive algebraic group action is equivariantly isomorphic to a closed $G$--subvariety of a finite dimensional complex vector space (see, for instance, Brion \cite[Proposition 1.9]{brion2010introduction}).
 \end{remark}

\subsection{Geometric invariant theory}\label{git}
Once a variety has been equipped with an algebraic group action, one might be tempted to take quotients.  Unfortunately, a naive construction of the quotient fails to produce a variety because the orbits of points are not always closed. Informally, geometric invariant theory remedies this situation by altering the quotient and gluing $``$bad orbits$"$ together. What follows is intended to be a brief reminder of this construction; the reader unfamiliar with these concepts is invited to consult Brion \cite{brion2010introduction} for the details.

Let $X$ be a $G$--variety and recall that this endows the ring of regular functions $\bb{C}[X]$ with a natural action of $G$. One would expect that the ring of regular functions of the quotient of $X$ by $G$ should correspond to those functions in $\bb{C}[X]$ which are invariant under $G$. 
Indeed, we define the affine Geometric Invariant Theory (GIT) quotient as the corresponding affine variety which we denote by $$
X/\!\!/G:=Specm(\bb{C}[X]^G).
$$
The topological space underlying $X/\!\!/G$ has a nice description of independent interest. The key observation here is that the closure of every $G$--\,orbit contains a unique closed orbit. Consequently, points of $X/\!\!/G$ correspond bijectively with the closed $G$--\,orbits. 
In fact, if we denote the set of closed $G$--\,orbits by $X/\!\!/_{top}\,G$ and consider the quotient map $\pi:X\raw X/\!\!/_{top}\,G$ sending $x\in X$ to the unique closed orbit in $\overline{G\cdot x}$, then $X/\!\!/_{top}\,G$ (equipped with the quotient topology) is homeomorphic to $X/\!\!/G$. 

\begin{remark}
This so-called topological Hilbert quotient $X/\!\!/_{top}\,G$ was investigated by Luna in \cite{luna1975certaines} and  \cite{luna1976fonctions} to understand the orbit space of real and complex algebraic group actions. In particular, he showed that $X/\!\!/_{top}\,G$ is the universal quotient in the category of Hausdorff spaces and in the category of complex analytic varieties.
\end{remark}

\begin{definition}
Let $X$ be a $G$--variety. A point $x\in X$ is said to be \emph{polystable} if its orbit $G\cdot x\subset X$ is Zariski closed. Since $G\cdot x$ is (Zariski) open in its (Zariski) closure, this is equivalent to requiring that $G\cdot x$ be closed in the Hausdorff topology. We denote the set of polystable points of $X$ by $X^{ps}$.
\end{definition}

\begin{example}
Let us consider the simple case where $G$  acts on itself by inner automorphisms. Since our definition of a $G$--variety requires an ambient vector space, we think of $G$ as a $G$--subvariety of the complex vector space of $n\times n$ matrices $G\subset M_n\bb{C}$ where $G$ acts on $M_n\bb{C}$ by conjugation. In this setting, it is well known that $x\in G$ is polystable if and only if $x$ is semisimple. 

Similarly, we can consider the diagonal action of $G$ on $G^r$ given by $g\cdot (x_1,\ldots, x_r):=(gx_1g\inv,\ldots,gx_rg\inv)$. Here, Richardson \cite[Theorem 3.6]{richardson1988conjugacy} generalized the previous characterization of polystability. The following theorem will play an important part in the proof of \hyperlink{Theorem III}{Theorem III} in Section \ref{kempfnessset}.
\begin{named}{Richardson's Theorem} \hypertarget{richpoly}{}
Let $G$ be a complex reductive linear algebraic group. For every tuple $(x_1,\ldots,x_n)\in G^n$, denote by $\mathscr{A}(x_1,\ldots,x_n)$ the algebraic subgroup of $G$ generated by the elements $\{x_1,\ldots,x_n\}$. If we consider the diagonal action of $G$ on $G^n$ by inner automorphisms, then $(x_1,\ldots,x_n)\in G^n$ is polystable if and only if $\mathscr{A}(x_1,\ldots,x_n)$ is reductive.
\end{named}
\end{example}
\subsection{The Kempf--Ness Theorem}\label{kntheorem}
Let us now consider the entire complex vector space $V$ as a $G$--variety via the algebraic inclusion  $G\subset \GL(V)$. If $K\subset G$ is a maximal compact subgroup, there is no loss of generality in assuming that $V$ is equipped with a $K$--\,invariant Hermitian inner product.
 It turns out that the associated norm $|\!|\cdot|\!|$ on $V$ sheds a lot of light on its polystable points. Indeed, one can determine if a vector $v\in V$ is polystable by understanding how its length changes as we move along its orbit $G\cdot v$. This variation in norm is captured by the \emph{Kempf--Ness function} which is defined for every  $v\in V$  by the rule \begin{equation}
\Psi_v:G\rightarrow \mathbb{R}\text{ , }g\mapsto |\!|g\cdot v|\!|^2.
\end{equation}
 Kempf and Ness  \cite[Theorem 0.2-0.3]{kempf1979length} used the convexity of these functions to formulate a simple characterization of polystability:

\begin{named}{Kempf--Ness Theorem}\hypertarget{KempfNess}{} 
Let $G\subset \GL(V)$ and suppose that $V$ is endowed with a $K$-invariant norm as above.
If $X\subset V$ is a $G$--subvariety, then:
\begin{enumerate}
 \item The point $x\in X$ is polystable if and only if $\Psi_x$ has a critical point.
 \item All critical points of $\Psi_x$ are minima.
 \end{enumerate}
 Moreover, if $\Psi_x(e)$ is a minimum then $K\cdot x=\{y\in G\cdot x:|\!|y|\!|=|\!|x|\!|\}$.
\end{named}
With this theorem in mind, we can state the most important definition of this paper.  As above, let $G\subset \GL(V)$ where $V$ is endowed with a $K$-invariant norm. Then, the \emph{Kempf--Ness set} of $V$ is defined as the $\mathcal{M}$inimal vectors 
\begin{equation}\mathcal{M}_V:=\{v\in V: |\!|g\cdot v|\!|\geq|\!|v|\!| \text{ for all }g\in G\}.\end{equation}
  Note that the \hyperlink{KempfNess}{Kempf--Ness Theorem} ensures that $\mathcal{M}_V\subset V^{ps}$ and, although $V^{ps}$ is an intrinsically defined subset of $V$, $\mathcal{M}_V$ depends on the choice of the invariant norm on $V$. If $X\subset V$ is a $G$--subvariety, we define the Kempf--Ness set of $X$ as $\mathcal{M}_X:=\mathcal{M}_V\cap X$.

The primary reason for calling the collection of minimal vectors the $``$Kempf--Ness set'' is the following corollary \cite[Corollary 4.7]{schwarz1989topology} of the \hyperlink{KempfNess}{Kempf--Ness Theorem}. It allows GIT quotients to be topologically represented via ordinary quotients and will play a r\^ole in the proof of \hyperlink{Theorem II}{Theorem II}.
\begin{corollary} \label{remarkable}
The geometric invariant theory quotient $X/\!\!/G$ is homeomorphic (in the Hausdorff topology) to the ordinary quotient $\mathcal{M}_X/K$.
\end{corollary} 

\begin{remark}
Whenever the Kempf--Ness set under consideration is clear from the context we shall omit the subscript and refer to it simply as $\mathcal{M}$.
\end{remark}

\subsection{The Neeman--Schwarz Theorem}\label{neemansection}

As above, let $G\subset \GL(V)$ act on the complex vector space $V$ and suppose the latter is endowed with a $K$--\,invariant norm $|\!|\cdot |\!|$. The following theorem \cite[Theorem 5.1]{schwarz1989topology} is the key known result used in the proof of our main theorems.
 \begin{named}{Neeman--Schwarz Theorem}\label{neemanSchwarztheorem}
\hypertarget{NeemanSchwarz}{} 
Let the complex vector space $V$ be a $G$--variety endowed with a $K$--\,invariant norm. If $\mathcal{M}$ denotes the Kempf--Ness set of $\,V$, then there is a $K$--\,equivariant strong deformation retraction of $\,V$ onto $\mathcal{M}$ which preserves $G$--orbits.

More precisely, there is  a map  $\varphi:V\times [0,1]\rightarrow V$
with the following properties:
\begin{enumerate}
\item $\varphi_0=\mathrm{Id}_V\text{, }\varphi_t|_\mathcal{M}=\mathrm{Id}_\mathcal{M}\text{, } \varphi_1(V)=\mathcal{M}$
\item $\varphi_t(v)\in G\cdot v$ for $0\leq t <1$
\item $\varphi_1(v)\in \overline{G\cdot v}$.
\end{enumerate}
\end{named}

For our purposes, the following corollary shall be most important:
\begin{corollary}
Let $X$ be a $G$--subvariety of $V$ and recall that $\mathcal{M}_X:=\mathcal{M}\cap X$. Restricting $\varphi$, we obtain a $K$--\,equivariant strong deformation retraction of $X$ onto $\mathcal{M}_X$. 
\end{corollary}

\begin{remark}
Initially, Neeman \cite[Theorem 2.1]{neeman1985topology} had proved  that there was a strong deformation retraction of $V$ onto a neighbourhood of $\mathcal{M}$. To do so, following a suggestion of Mumford, he considered the negative gradient flow of a function assigning to each $v\in V$ the $``$size'' of the differential $(d\Psi_v)_e$.
However, he believed that the local flow could be extended to a retraction of $V$ onto $\mathcal{M}$ and conjectured a functional inequality that would allow him to do so.  Later, he realized that this was in fact a special case of an inequality due to  Lojasiewicz \cite{lojasiewicz1965ensembles}. Finally a complete account of the proof of the \hyperlink{NeemanSchwarz}{Neeman--Schwarz Theorem} was given by Schwarz in  \cite[Theorem 5.1]{schwarz1989topology}. 

\end{remark}
\begin{remark}
Although we have been working over $\bb{C}$, Kempf--Ness theory was extended to the real setting by Richardson and Slodowy. For example, in the case of  \hyperlink{richpoly}{Richardson's Theorem} it is shown in Richardson \cite[Theorem 11.4]{richardson1988conjugacy}   that 
$(x_1,\ldots,x_n)\in G(\kern 0.05 em \bb{R})^n$ is polystable for the action of $G(\mathbb{R})$ if and only if $(x_1,\ldots,x_n)\in G^n$ is polystable for the action of $G$.
For the \hyperlink{KempfNess}{Kempf--Ness Theorem}, its Corollary and the \hyperlink{NeemanSchwarz}{Neeman--Schwarz Theorem}, it is shown in Richardson--Slodowy \cite{richardson1990minimum} that the same results hold when one replaces $G$ by its group of real points $G(\kern 0.05 em \bb{R})$. This will be elaborated upon in Section \ref{real} when we prove the real case of \hyperlink{Theorem I}{Theorem I}.
\end{remark}

\section{The representation variety}\label{the representation variety}
 Let $\Gamma$ be a finitely generated group and let $G$ be a complex reductive linear algebraic group. In a coarse sense, this paper is concerned with the finite dimensional representation theory of $\Gamma$ in $G$. When $\Gamma$ is finite, this essentially reduces to the correspondence between linear representations and characters due to Frobenius. When $\Gamma$ is infinite, the analogue of character theory is the parametrization of representations by geometric varieties. In this section, we study the topology of such varieties.
 
In order to endow the space   $\Hom(\kern 0.05em\Gamma,G\kern 0.05em)$ of homomorphisms $\rho\negthinspace:\Gamma\raw G$ with a geometric structure, it is convenient to work with a fixed generating set $\gamma_1,\ldots,\gamma_r$ of $\Gamma$.  Since any representation $\rho\negthinspace:\Gamma\raw G$ is uniquely determined by the image $\rho(\gamma_i)$ of its generators, the evaluation map 
\begin{equation}\Hom(\kern 0.05em\Gamma,G\kern 0.05em)\raw G\times\ldots\times G=G^r\text{, }\rho\mapsto (\rho(\gamma_1),\ldots,\rho(\gamma_r)) \label{evaluation}\end{equation} is a bijection between $\Hom(\kern 0.05em\Gamma,G\kern 0.05em)$ and a Zariski closed subset of $G^r$ carved out by the relations of $\Gamma$. It is well known (see, for instance, Lubotzky and Magid \cite{lubotzky1985varieties}) and otherwise easy to see that the geometric structure induced on $\Hom(\kern 0.05em\Gamma,G\kern 0.05em)$ by this bijection doesn't depend on the chosen presentation of $\Gamma$.

 \begin{remark}
With the hope that no confusion shall arise, we denote elements of $\Hom(\kern 0.05em\Gamma,G\kern 0.05em)$ by $\rho$ when we consider them as homomorphisms and by tuples $(m_1,\ldots,m_r)$ when we identify them with points in $G^r$.
\end{remark}

The action of $G$ on itself by conjugation induces a diagonal action of $G$ on the variety $\Hom(\kern 0.05em\Gamma,G\kern 0.05em)\subset G^r$. Explicitly, this action is given by the rule $$G\times \Hom(\kern 0.05em\Gamma,G\kern 0.05em)\raw\Hom(\kern 0.05em\Gamma,G\kern 0.05em)\text{, }g\cdot (m_1,\ldots,m_r):=(gm_1g\inv,\ldots,gm_rg\inv).$$
 Since two representations are in the same $G$--orbit when their image differs by a change of basis, one is often interested in classifying them modulo the action of $G$. This can be done by considering the geometric invariant theory quotient $\Hom(\kern 0.05em\Gamma,G\kern 0.05em)/\!\!/G$ which happens to parametrize isomorphism classes of completely reducible representations. Describing these varieties geometrically might be viewed as analogous to determining the  characters of a finite group.

As an affine variety, $\Hom(\kern 0.05em\Gamma,G\kern 0.05em)$ also inherits a natural Hausdorff topology given from an embedding into affine space. In fact, it coincides with the subspace topology of $\Hom(\kern 0.05em\Gamma,G\kern 0.05em)\subset G^r$ when $G$ is viewed as a Lie group. If $K\subset G$ is any subgroup, there is a natural inclusion $\Hom(\kern 0.05em\Gamma,K\kern 0.05em)\hookrightarrow\Hom(\kern 0.05em\Gamma,G\kern 0.05em)$ so we endow $\Hom(\kern 0.05em\Gamma,K\kern 0.05em)$ with the subspace topology. 
  \hyperlink{Theorem I}{Theorem I} states that under certain assumptions  there is a deformation retraction of $\Hom(\kern 0.05em\Gamma,G\kern 0.05em)$ onto $\Hom(\kern 0.05em\Gamma,K\kern 0.05em)$. The goal of this section is to show how the algebro-geometric tools developed in the previous section can be used in this setting.

\subsection{The Frobenius norm}

\label{frobeniusnorm}
Let $G$ be a (possibly disconnected) complex reductive linear algebraic group and let $K\subset G$ be a maximal compact subgroup. In this subsection, we complete an intermediate step in the proof of \hyperlink{Theorem I}{Theorem I} by showing how the \hyperlink{NeemanSchwarz}{Neeman--Schwarz Theorem} can be applied to $\Hom(\kern 0.05em\Gamma,G\kern 0.05em)$.  Keeping with the general spirit of this paper, instead of working abstractly with our algebraic group $G$, we 
 choose a matrix representation $G\hookrightarrow \sln$ for which $K=G\cap \SU_n$.

In order to apply the algebro-geometric results from Section \ref{algebraic actions}, we need to endow $\Hom(\kern 0.05em\Gamma,G\kern 0.05em)$ with the structure of an affine $G$--variety.  To do this, consider $G\subset \sln$ as a variety in the vector space  $M_n\bb{C}$ of $n\times n$ complex matrices.  Using the natural identifications arising from the evaluation map $(\ref{evaluation})$, we can realize $\Hom(\kern 0.05em\Gamma,G\kern 0.05em)$ as an embedded variety  
$$\Hom(\kern 0.05em\Gamma,G\kern 0.05em)\subset G^r\subset \oplus_{j=1}^r M_n\bb{C}.$$
Once this has been done, the diagonal action of $G\subset \sln$ on $\oplus_{j=1}^r M_n\bb{C}$ by simultaneous conjugation  endows $\Hom(\kern 0.05em\Gamma,G\kern 0.05em)$ with the structure of a $G$--subvariety of the vector space $\oplus_{j=1}^r M_n\bb{C}\cong \bb{C}^{rn^2}$.

In this setting, there is a prototypical $K$--\,invariant norm on $\oplus_{j=1}^r M_n\bb{C}$.
  Recall that the vector space $M_n\bb{C}$ admits a natural Hermitian inner product called the \emph{Frobenius inner product} (see, for instance, Horn and Johnson \cite[Chapter 5]{horn1990matrix}). If we denote the adjoint of a matrix $a$ by $a^\ast:=\overline{a}^{\,t}$, then this inner product is defined by $\langle a,b\rangle:=\mathrm{trace}(a^\ast b)$ for $a,b\in M_n\bb{C}$. It is easy to check that for any unitary matrix $u\in \U_n$ we have $\langle uau\inv,ubu\inv\rangle=\langle a,b\rangle$.  The Frobenius inner product can then be extended $``$coordinate wise'' to $\oplus_{j=1}^r M_n\bb{C}$ by the rule $$\langle (a_1,\ldots, a_r),(b_1,\ldots,b_r)\rangle:=\sum_{j=1}^r\langle a_j,b_j\rangle$$
 and the associated norm on $\oplus_{j=1}^r M_n\bb{C}$ is unitary invariant (and therefore $G\,\cap\, \SU_n=K$--\,invariant).

\begin{remark} We henceforth refer to the $K$--\,invariant norm described above as the Frobenius norm of the $G$--variety $\Hom(\kern 0.05em\Gamma,G\kern 0.05em)\subset \oplus_{j=1}^r M_n\bb{C}$.\end{remark} 

 Once it has been equipped with this additional structure, we can define the Kempf--Ness set $\mathcal{M}$ of $\Hom(\kern 0.05em\Gamma,G\kern 0.05em)$ as in Section \ref{neemansection} to be those representations that are of minimal norm within their $G$--\,orbit.
 The \hyperlink{NeemanSchwarz}{Neeman--Schwarz Theorem} now provides us with a $K$--\,equivariant deformation retraction of $\Hom(\kern 0.05em\Gamma,G\kern 0.05em)$ onto  $\mathcal{M}$.

\subsection{The Kempf--Ness set}\label{the Kempf--Ness set}
In order to fruitfully apply the \hyperlink{NeemanSchwarz}{Neeman--Schwarz Theorem} to $\Hom(\kern 0.05em\Gamma,G\kern 0.05em)$, we need to understand the target of the retraction it provides. For a general group $\Gamma$, e.g., if $\Gamma$ is a free group, explicitly determining the Kempf--Ness set of $\Hom(\kern 0.05em\Gamma,G\kern 0.05em)$ appears to be a hopeless task. This is why we henceforth let $\Gamma$ be a finitely generated \emph{nilpotent} group. In this case, the Kempf--Ness set admits a very nice description:
\begin{named}{\hyperlink{Theorem III}{Theorem III}}\label{kempfnessset}
Let $G\subset \sln$ be a complex reductive algebraic group and suppose that $K=G\cap \SU_n\subset G$ is a maximal compact subgroup. If $\Gamma$ is a finitely generated nilpotent group and $\mathcal{M}$ denotes the Kempf--Ness set of  the $G$--subvariety $\Hom(\kern 0.05em\Gamma,G\kern 0.05em)\subset \oplus_{j=1}^r M_n\bb{C}$ equipped with the Frobenius norm, then $$
\mathcal{M}=\{\rho\in\Hom(\kern 0.05em\Gamma,G\kern 0.05em):\rho(\Gamma)\text{ consists of normal matrices}\}.
$$
\end{named}

The proof of \hyperlink{Theorem III}{Theorem III} relies on two lemmas. Before diving into the heart of the matter, we refresh the reader's memory with some elementary linear algebra.  
Recall that a matrix $m\in M_n\bb{C}$ is \emph{normal} if $mm^\ast=m^\ast m$ or, equivalently, $V$ has an orthogonal basis of eigenvectors of $m$.  The following well known\footnote{We invite the intrigued reader to consult Grone--Johnson--Sa--Wolkowicz  \cite{grone1987normal} where this is the $53$rd (out of $70$!) characterization of normality.} characterization of normal matrices follows from the unitary invariance of the Frobenius norm $|\!|\cdot|\!|$ of $M_n\bb{C}$ and the Schur Triangularization Theorem.

\begin{lemma}\label{morse}
If $\lambda_1,\ldots,\lambda_n$ denote the eigenvalues of the matrix $m\in M_n\bb{C}$, then  $$|\!|m|\!|^2\geq \sum_{j=1}^n |\lambda_j|^2.$$ Moreover, equality holds if and only if $m$ is a normal matrix.\qed
\end{lemma}

Our second lemma  concerns the images of certain  representations of $\Gamma$ in $G$:

\begin{lemma}\label{oldbadlemma}
Let $G\subset \sln$ be a reductive algebraic group and suppose that $K=G\cap \SU_n$ is a maximal compact subgroup of $G$. If
$N\subset G$ is a nilpotent reductive algebraic subgroup, then there is some $g\in G$ for which $gNg\inv$ consists of normal matrices.
\end{lemma}
\begin{proof}
Recall from Lemma \ref{componentnilp} that a reductive nilpotent algebraic group $N$ has a unique maximal compact subgroup $C$. Since $C$ is contained in a maximal compact subgroup of $G$ and all such groups are conjugate in $G$, there is some $g\in G$ for which $gCg\inv\subset K$. We can therefore assume (after conjugating by $g$) that $N\subset G\subset\sln$ and $C\subset K\subset \SU_n$.

Being a reductive algebraic group, $N$ is the complexification of $C$. Thus, it admits a polar decomposition (c.f. Section \ref{complexreductive}):$$N=C\cdot exp(i\fk{c})\simeq C\times exp(i\fk{c})$$ where $exp(i\fk{c})\subset N^o$ is central in $N$ by Lemma \ref{componentnilp}. This shows that any matrix $m\in N$ can be written uniquely as a product $m=k\cdot h$ for some unitary $k\in C$ and Hermitian $h\in exp(i\fk{c})$. Since unitary and Hermitian matrices are normal and since the product of two commuting normal matrices is normal, $N$  consists of normal matrices.
\end{proof}
The first lemma is used to show that normal representations are critical points of the Kempf--Ness functions so they belong to the Kempf--Ness set. Then, the second lemma shows that polystable representations become normal after conjugation so the Kempf--Ness set consists exclusively of normal representations. More precisely:
\begin{proof}[Proof of \hyperlink{Theorem III}{Theorem III}]
Let $\rho\in\Hom(\kern 0.05em\Gamma,G\kern 0.05em)$ and consider the corresponding tuple $$(m_1,\ldots,m_r):=(\rho(\gamma_1),\ldots,\rho(\gamma_r))\subset G^r.$$ If $m_j$ is a normal matrix for each $1\leq j\leq r$, then Lemma \ref{morse} ensures that $$|\!|(m_1,\ldots,m_r)|\!|\leq |\!|(gm_1g\inv,\ldots,gm_rg\inv)|\!|$$ for any $g\in G$. Consequently, the representation $\rho$ is in $\mathcal{M}$; this shows that $$\mathcal{M}\supset\{\rho\negthinspace:\rho(\gamma_j)\text{ is normal for }1\leq j\leq r\}\supset\{\rho\negthinspace:\rho(\Gamma)\text{ consists of normal matrices}\}.$$
We now show that the three sets coincide. 
Seeking a contradiction, suppose that $\rho\in \mathcal{M}\subset \Hom(\kern 0.05em\Gamma,G\kern 0.05em)^{\,ps}$ is a representation whose image doesn't consist of normal matrices.  Recall that \hyperlink{richpoly}{Richardson's Theorem} characterizes polystable points in $G^r$ for the conjugation action of $G$ as those tuples $(m_1,\ldots, m_r)$ for which $\mathscr{A}(m_1,\ldots,m_r)\subset G$ is reductive. Therefore, the Zariski closure of the image of our representation $$N:=\mathscr{A}(\rho(\Gamma))=\mathscr{A}(\rho(\gamma_1),\ldots,\rho(\gamma_r))$$ is a reductive nilpotent algebraic subgroup of $G$. But then, Lemma \ref{oldbadlemma} produces an element $g\in G$ for which $gNg\inv$ consists of normal matrices and, by Lemma \ref{morse}, this contradicts the assumption that $\rho$ was of minimal norm within its $G$--\,orbit.
\end{proof}

\begin{example}\label{counter}
While it is more or less clear that if we replace $\Gamma$ by a free group then $\mathcal{M}$ does not consist of normal representations (see, for instance, Florentino--Lawton \cite[Remark 3.2]{florentino2014topology}), it may be somewhat surprising that this already fails to be the case when $\Gamma$ is solvable.   To see this, consider the solvable group $\Gamma:=\bb{Z}\rtimes \bb{Z}/4\bb{Z}$ generated by the elements $(1,0)$, $(0,1)$ and $(1,1)$.  For this generating set, an obvious  representation $\rho\in\Hom(\kern 0.05 em\Gamma,\SL_2\!\bb{C}\kern 0.05 em)$ is given by mapping 
\[
(1,0)\mapsto \left(\begin{array}{cc}\lambda & 0 \\0 & \lambda\inv\end{array}\right),\,(0,1)\mapsto   \left(\begin{array}{cc}0 & 1 \\-1 & 0\end{array}\right)\text{ and }(1,1)\mapsto \left(\begin{array}{cc}0 & \lambda \\-\lambda\inv & 0\end{array}\right)
\]
 where $|\lambda|\neq 1,0$. A simple computation shows that $\rho$ lies in the Kempf--Ness set while the matrix $\rho(1,1)$  fails to be normal. This example was shown to us by Lior Silberman.

\end{example}

\section{The scaling operation}\label{the scaling operation}
In the previous section, we determined the Kempf--Ness set of the representation variety $\Hom(\kern 0.05em\Gamma,G\kern 0.05em)$. To complete the proof of \hyperlink{Theorem I}{Theorem I} in the complex case, we seek a deformation retraction of this set onto $\Hom(\kern 0.05em\Gamma,K\kern 0.05em)$. The goal of this section is to show that scaling eigenvalues does the trick.
\subsection{Scaling eigenvalues}
\label{scalingeigenvalues}
Following Pettet and Souto \cite{pettet2011commuting}, we define an eigenvalue scaling operator on the semisimple elements of an algebraic group. Let $\bb{D}_n$ denote the diagonal subgroup of $\sln$, identified in the usual way with a subgroup of $(\bb{C}^\times)^n$. Consider the map $$\tilde{\sigma}:[0,1]\times \bb{D}_n\raw \bb{D}_n\text{, }(t,(\lambda_i))\mapsto \tilde{\sigma}_t((\lambda_i)):=(e^{-t\log|\lambda_i|}\lambda_i)$$
where $\log(\cdot)$ denotes the standard real logarithm. Notice that $\tilde{\sigma}$ is continuous, equivariant under the action of the Weyl group, and that $
\tilde{\sigma}_0(g)=g\text{ while }\tilde{\sigma}_1(g)\in\bb{D}_n\cap \SU_n
$
 for all $g\in\bb{D}_n$.
 
 Given now a semisimple element $g\in\sln$, there is some  $h\in \sln$ for which $hgh\inv\in\bb{D}_n$. Since $\tilde{\sigma}$ is invariant under the Weyl group, we see that for all $t\in[0,1]$ the element $$\sigma_t(g):=h\inv\tilde{\sigma}_t(hgh\inv)h$$ is independent of the choice of $h$. In particular, we obtain a well defined \emph{scaling map} \begin{equation}
 \sigma:[0,1]\times (\sln)_s\raw (\sln)_s\text{, }(t,g)\mapsto\sigma_t(g)
 \end{equation} where $(\sln)_s$ denotes the set of semisimple elements in $\sln$. This map preserves commutativity and scales the eigenvalues of the semisimple elements of $\sln$ until they land in a subgroup conjugate to $\SU_n$. Using this map, the following lemma is proved in \cite[Section 8.3]{pettet2011commuting}:
 
\begin{lemma}[Pettet--Souto]\label{scalesouto}
If $G$ is any linear algebraic group and $K\subset G$ is a maximal compact subgroup,  then there is a $G$--equivariant continuous map
$$\sigma:[0,1]\times G_s\raw G_s\text{, }(t,g)\mapsto \sigma_t(g) $$
which satisfies the following properties for all $g,g_1,g_2\in G_s$:\begin{enumerate}
\item $\sigma_0(g)=g$ and $\sigma_1(g)$ lies in a compact subgroup of $G$ conjugate to $K$.
\item $\forall$ $ 0\leq t\leq 1$, $\sigma_t$ fixes pointwise the subgroups of $G$ conjugate to $K$.
\item $\forall$ $ 0\leq t\leq 1$, if $g_1g_2=g_2g_1$  then $\sigma_t(g_1g_2)=\sigma_t(g_1)\sigma_t(g_2)=\sigma_t(g_2)\sigma_t(g_1).$

\item If $G$ is defined over $\bb{R}$ then $\sigma$ is equivariant under complex conjugation.
\end{enumerate}
\end{lemma}

\subsection{Scaling the Kempf--Ness set}
Let us now return to the setting of \hyperlink{Theorem I}{Theorem I} where $\Gamma$ is a finitely generated nilpotent group and $G$ is a complex reductive algebraic group.
In Section \ref{frobeniusnorm}, given a maximal compact subgroup $K\subset G$, we chose a faithful representation $G\hookrightarrow \SL_n\!\bb{C}$ for which $K=G\cap \SU_n$ to endow $\Hom(\kern 0.05em\Gamma,G\kern 0.05em)\subset \oplus_{j=1}^r M_n\bb{C}$ with a $K$--\,invariant Frobenius norm. This allowed us to define the Kempf--Ness set $\mathcal{M}$ of $\Hom(\kern 0.05em\Gamma,G\kern 0.05em)$ as the set of representations of minimal norm within their $G$--\,orbit. In fact, \hyperlink{Theorem III}{Theorem III} showed that 
$$
\mathcal{M}=\{\rho\in\Hom(\kern 0.05em\Gamma,G\kern 0.05em):\rho(\Gamma)\text{ consists of normal matrices}\}.
$$
The goal of this subsection is to use the scaling map $\sigma$ from Lemma \ref{scalesouto} to prove the following proposition:  

\begin{proposition}\label{nilpotentscaling}
Keeping the notation as above, there is a $K$--\,equivariant strong deformation retraction of the Kempf--Ness set $\mathcal{M}$ onto $\Hom(\kern 0.05em\Gamma,K\kern 0.05em)$.
\end{proposition}
Before proving the proposition, we need to understand how the image of representations in the Kempf--Ness set behave under scaling; this is the content of the following lemma.

\begin{lemma}\label{scalingrepresentations}
Let $N\subset \sln$ be a nilpotent reductive algebraic group. If $\sigma_t$ denotes the scaling map from \emph{Lemma \ref{scalesouto}}, then the map $N\mapsto \sigma_t(N)$ is a group homomorphism.
\end{lemma}
\begin{proof}
Recall from Lemma \ref{componentnilp} that a nilpotent reductive algebraic group $N$  consists of semisimple elements and has a unique maximal compact subgroup $C$. Moreover, it admits a polar decomposition  
$$N=C\cdot exp(i\fk{c})\simeq C\times exp(i\fk{c})$$ where $\fk{c}$ denotes the Lie algebra of $C$ and $exp(i\fk{c})\subset N^o$ is central in $N$ by Lemma \ref{componentnilp} again.

 Let $m_1=k_1\cdot h_1$ and $m_2=k_2\cdot h_2$ denote the unique decomposition of elements $m_1,m_2\in N$ into products of elements $k_1,k_2\in C$ and $h_1,h_2\in exp(i\fk{c})$.  Notice that since $C$ is maximal compact it is fixed pointwise by $\sigma_t$. Then, since $h_1$ and $h_2$ are central in $N$,  we have:
$$\sigma_t(m_1 m_2)=\sigma_t(k_1 h_1 k_2 h_2)=\sigma_t(k_1 k_2 h_1 h_2)=\sigma_t(k_1 k_2)\sigma_t(h_1) \sigma_t(h_2)=k_1k_2\sigma_t(h_1)\sigma_t(h_2)$$ 
$$=k_1\sigma_t(h_1)k_2\sigma_t(h_2)=\sigma_t(k_1)\sigma_t(h_1)\sigma_t(k_2)\sigma_t(h_2)=\sigma_t(k_1h_1)\sigma(k_2h_2)=\sigma_t(m_1)\sigma_t(m_2).
$$
Here, we repeatedly use parts $(2)$ and $(3)$ of Lemma \ref{scalesouto}.
\end{proof}
\begin{proof}[Proof of Proposition \ref{nilpotentscaling}]
Let us denote by $\mathcal{N}$ the set of normal matrices in $\sln$ and consider the restriction of $\sigma_t$ to $\mathcal{N}$. Observe that $\sigma_0(\mathcal{N})=\mathcal{N}$ while $\sigma_1(\mathcal{N})=\SU_n$ because a matrix is unitary if and only if it is normal and its eigenvalues all have unit norm. Since the representations $\rho$ in the Kempf--Ness set $\mathcal{M}$ of $\Hom(\kern 0.05em\Gamma,G\kern 0.05em)$ are precisely those for which $N:=\mathscr{A}(\rho(\Gamma))\subset \mathcal{N}$ is a reductive nilpotent group, it follows from Lemma \ref{scalingrepresentations} that $\sigma_t\circ \rho$ is also in $\mathcal{M}$. Consequently, the map $\psi:\mathcal{M}\times [0,1]\raw \mathcal{M}$ where $\psi(\rho,t):=\sigma_t\circ \rho$ is a strong deformation retraction with $\psi(\rho,1)(\Gamma)\subset K=G\cap \SU_n$. 
\end{proof}

\section{The real and complex cases}\label{real}
The previous two sections contain a complete proof of \hyperlink{Theorem I}{Theorem I} when $G$ is a complex reductive linear algebraic group. To recapitulate, once $\Hom(\kern 0.05em\Gamma,G\kern 0.05em)$ has been endowed with the appropriate structure, the retraction proceeds essentially in two steps:
$$
\Hom(\kern 0.05em\Gamma,G\kern 0.05em)\xrightarrow[Theorem]{Neeman-Schwarz}\text{ {Kempf--Ness Set} }\xrightarrow[Eigenvalues]{Scaling}\Hom(\kern 0.05em\Gamma,K\kern 0.05em).
$$ 
In this section, we use real Kempf--Ness theory to adapt these  steps when $G$ is replaced by its group of real points $G(\kern 0.05 em \bb{R})$. 

\subsection{Real Kempf--Ness theory}
To begin, we discuss the work of Richardson and Slodowy found in    \cite{richardson1990minimum} where they develop real analogues of the main results in Kempf--Ness theory.

Let $V$ be a finite dimensional complex vector space with real structure $V(\kern 0.05 em \bb{R})$ and let $H$ be a positive-definite Hermitian form on $V$. In this case, we always have $H(x,y)=S(x,y)+iA(x,y)$ where $S$ (respectively $A$) is a symmetric (respectively alternating) real-valued $\bb{R}$-bilinear form on $V$. We say that the Hermitian form $H$ is \emph{compatible} with the $\bb{R}-$structure of $V$ if $A$ vanishes on $V(\kern 0.05 em\bb{R})\times V(\kern 0.05 em \bb{R})$.

\begin{example}
The Frobenius inner product on $M_n\bb{C}$  is \emph{compatible} with the $\bb{R}-$structure $M_n\bb{R}$.
\end{example}

Let us suppose now that the complex reductive algebraic group $G\subset \GL(V)$ is defined over $\bb{R}$ and denote its group of real points by $G(\kern 0.05 em \bb{R})$.  Since $G$ is reductive, it follows from the work of Mostow \cite{mostow1955self} that, upon identifying $V\cong \bb{C}^m$ and $\GL(V)\cong \GL_m\bb{C}$, we can assume $G$ to be stable under the involution $\theta(g):=(g^\ast)\inv$. In this case, the fixed points $$K:=G^\theta=\{g\in G:\theta(g)=g\}$$ form a maximal compact subgroup of $G$ and the real fixed points $
K(\kern 0.05 em \bb{R})=G(\kern 0.05 em \bb{R})^\theta
$ form a maximal compact subgroup of $G(\kern 0.05 em \bb{R})$. One can show \cite[Appendix 2]{richardson1990minimum} that the vector space $V$ may always be equipped with a $K$--\,invariant Hermitian inner product which is compatible with the $\bb{R}-$structure $V(\kern 0.05 em \bb{R})$. We denote the associated $K$--\,invariant norm on $V$ by $|\!|\cdot |\!|$. 

The (complex) Kempf--Ness set $\mathcal{M}$ of $V$ for the action of $G$ can now be defined as in Section \ref{kntheorem} 
and the (real) Kempf--Ness set $\mathcal{M}(\kern 0.05 em \bb{R})$ of $V(\kern 0.05 em \bb{R})$ for the action of $G(\kern 0.05 em \bb{R})$ can be defined by the analogous rule  
\begin{equation}\mathcal{M}(\kern 0.05 em \bb{R}):=\{v\in V(\kern 0.05 em \bb{R}): |\!|g\cdot v|\!|\geq|\!|v|\!| \text{ for all }g\in G(\kern 0.05 em \bb{R})\}.\end{equation}
The following lemma is an immediate consequence of \cite[Lemma $8.1$]{richardson1990minimum}:
\begin{lemma}\label{intersectionkempfness}
$\mathcal{M}(\kern 0.05 em \bb{R})=\mathcal{M}\cap V(\kern 0.05 em \bb{R}).$
\end{lemma}
At this point we can state the two only results of Richardson and Slodowy that we shall need. The first is an analogue of Corollary \ref{remarkable}   \cite[Theorem $7.7$]{richardson1990minimum}:
\begin{theorem}\label{realkempf}
Keeping the notation as above, let $X$ be a closed $G(\kern 0.05 em \bb{R})-$stable subset of $V(\kern 0.05 em \bb{R})$. If $\mathcal{M}_X(\kern 0.05 em \bb{R}):=\mathcal{M}(\kern 0.05 em \bb{R})\cap X$, then there is a homeomorphism $\mathcal{M}_X(\kern 0.05 em \bb{R})/K(\kern 0.05 em \bb{R})\cong X/\!\!/_{top}\,G(\kern 0.05 em \bb{R})$.
\end{theorem}
\begin{remark}
Here, $X/\!\!/_{top}\,G(\kern 0.05 em \bb{R})$ denotes the \emph{topological Hilbert quotient} (c.f. Section \ref{git} and  \cite[Section 7.2]{richardson1990minimum}).
\end{remark}
The second is an analogue of the \hyperlink{NeemanSchwarz}{Neeman--Schwarz Theorem} \cite[Theorem $9.1$]{richardson1990minimum}:
\begin{theorem}\label{realneeman}
Keeping the notation as above, there is a continuous $K(\kern 0.05 em \bb{R})-$equivariant strong deformation retraction $\varphi:V(\kern 0.05 em \bb{R})\times [0,1]\raw V(\kern 0.05 em \bb{R})$ of $V(\kern 0.05 em \bb{R})$ onto $\mathcal{M}(\kern 0.05 em \bb{R})$. Moreover, the deformation is along orbits of $G(\kern 0.05 em \bb{R})$, that is, $\varphi_t(v)\subset G(\kern 0.05 em \bb{R})\cdot v$ for $0\leq t<1$ and $\varphi_1(v)\in \overline{G(\kern 0.05 em \bb{R})\cdot v}$.
\end{theorem}
\begin{remark}
As in the complex case, if $X$ is a closed $G(\kern 0.05 em \bb{R})-$stable subset of $V(\kern 0.05 em \bb{R})$, then the restriction of $\varphi$ to $X\times[0,1]$ is a deformation retraction of $X$ onto $\mathcal{M}_X(\kern 0.05 em \bb{R}):=\mathcal{M}(\kern 0.05 em \bb{R})\cap X$.
\end{remark}

\subsection{Proofs of {Theorem I} and {Theorem II}}

We are now ready to complete the proofs of our main theorems. We begin with the proof of the real case of \hyperlink{Theorem I}{Theorem I}:
\begin{named}{\hyperlink{Theorem I}{Theorem I}}
{Let $\Gamma$ be a finitely generated nilpotent group and let $G$ be the group of complex or real points of a reductive linear algebraic group defined over $\bb{R}$. If $K$ is a maximal compact subgroup of $G$, then there is a $K$--\,equivariant strong deformation retraction of $\Hom(\kern 0.05em\Gamma,G\kern 0.05em)$ onto $\Hom(\kern 0.05em\Gamma,K\kern 0.05em)$.}
\end{named}
\begin{proof}[Proof of the real case] 
Let $G\subset\sln$ be a complex reductive algebraic group defined over $\bb{R}$ and denote its group of $\bb{R}-$points by $G(\kern 0.05 em \bb{R})$.    In view of the Conjugacy Theorem \ref{conjugacy}, it is no loss of generality to assume that $K=G\cap \SU_n$ is a maximal compact subgroup of $G$ for which  $K(\kern 0.05 em \bb{R})=K\cap G(\kern 0.05 em \bb{R})$ is the maximal compact subgroup of $G(\kern 0.05 em \bb{R})$ in the statement of \hyperlink{Theorem I}{Theorem I}.
We can then proceed essentially as in Section \ref{frobeniusnorm}. If $\Gamma$ is generated by $r$ elements, we embed $\Hom(\kern 0.05 em\Gamma,G(\kern 0.05 em \bb{R})\kern 0.05 em)\subset \oplus_{j=1}^rM_n\bb{R}$ and Theorem \ref{realneeman} provides us with a $K(\kern 0.05 em \bb{R})-$equivariant strong deformation retraction of $\Hom(\kern 0.05 em\Gamma,G(\kern 0.05 em \bb{R})\kern 0.05 em)$ onto its Kempf--Ness set $\mathcal{M}(\kern 0.05 em \bb{R})$. Here, Lemma \ref{intersectionkempfness} and \hyperlink{Theorem III}{Theorem III} imply that $$\mathcal{M}(\kern 0.05 em \bb{R})=\{\rho\in\Hom(\kern 0.05 em\Gamma,G(\kern 0.05 em \bb{R})\kern 0.05 em):\rho(\Gamma)\text{ consists of normal matrices }\}.$$
The proof is then completed by Lemma \ref{scalesouto} and Proposition \ref{nilpotentscaling} which show that scaling eigenvalues induces a  $K(\kern 0.05 em \bb{R})-$equivariant retraction of $\mathcal{M}(\kern 0.05 em \bb{R})$ onto $\Hom(\kern 0.05 em\Gamma,K(\kern 0.05 em \bb{R})\kern 0.05 em)$.
\end{proof}

Finally, we prove the analogous theorem for character varieties: 
\begin{named}{\hyperlink{Theorem II}{Theorem II}}
Let $\Gamma$ be a finitely generated nilpotent group and let $G$ be the group of complex or real points of a reductive linear algebraic group, defined over $\bb{R}$ in the latter case. If $K\subset G$ is any maximal compact subgroup, then there is a strong deformation retraction of the character variety $\Hom(\kern 0.05em\Gamma,G\kern 0.05em)/\!\!/G$ onto $\Hom(\kern 0.05em\Gamma,K\kern 0.05em)/K$.
\end{named}

\begin{proof}
Let $G\subset\SL_n\!\bb{C}$ be a complex reductive linear algebraic group and equip $\Hom(\kern 0.05em\Gamma,G\kern 0.05em)\subset G^r$ with the Frobenius norm as in Section \ref{frobeniusnorm}. If we let $\mathcal{M}$ denote its Kempf--Ness set, then   
Corollary \ref{remarkable}  shows that there is a homeomorphism
$\Hom(\kern 0.05em\Gamma,G\kern 0.05em)/\!\!/G\cong \mathcal{M}/K$. Proposition \ref{nilpotentscaling} then produces a $K$--equivariant strong deformation retraction of $\mathcal{M}$ onto $\Hom(\kern 0.05em\Gamma,K\kern 0.05em)$ which induces a strong deformation retraction of $\mathcal{M}/K$ onto  $\Hom(\kern 0.05em\Gamma,K\kern 0.05em)/K$.
If $G$ is defined over $\bb{R}$, the argument above can be carried out for its group of real points $G(\kern 0.05 em \bb{R})$ by viewing $\Hom(\kern 0.05 em\Gamma,G(\kern 0.05 em \bb{R})\kern 0.05 em)/\!\!/_{top} \,G(\kern 0.05 em \bb{R})$ as a topological Hilbert quotient (c.f. Section \ref{git} and Richardson--Slodowy \cite[Section 7.2]{richardson1990minimum}) and invoking Theorem \ref{realkempf} instead of  Corollary \ref{remarkable}. 
\end{proof} 

\appendix
\section{Expanding nilpotent groups}\label{appendix}
In this appendix, we complete the comment made in the introduction on the necessity of the assumption that $G$ be reductive in \hyperlink{Theorem I}{Theorem I}.  The following content was shown to us by  Lior Silberman and Juan Souto. We refer the reader to Raghunathan  \cite[Chapter 2]{raghunathan1972discrete} for generalities on nilpotent lattices in unipotent Lie groups.

 Recall from the work of Mal'cev that every finitely generated torsion free nilpotent group admits a canonical completion to a unipotent Lie group which is usually called its \emph{Mal'cev completion}. More precisely, one has the following theorem:

\begin{theorem}[Mal'cev]
Let $\Gamma$ be a finitely generated torsion free nilpotent group. Then, there is a canonical unipotent Lie group $U$ for which $\Gamma\subset U$ is a lattice. Moreover, if $U'$ is any unipotent Lie group, then any homomorphism $\Gamma\raw U'$ extends uniquely to a homomorphism $U\raw U'$.
\end{theorem}

Our goal is to show that one can choose $\Gamma$ and $U$ in a way that makes the conclusion of \hyperlink{Theorem I}{Theorem I} fail when we consider the space of representations $\Hom(\kern 0.05em\Gamma,U\kern 0.05em)$. It turns out that the key property of $\Gamma\subset U$ in this setting is the $``$richness$"$ of its automorphism group as described in the following definition.

\begin{definition}[Definition.]
A finitely generated torsion free nilpotent group $\Gamma$ is said to be \emph{expanding} if the Lie algebra $\mathfrak{u}$ of its Mal'cev completion $U$ admits a semisimple automorphism, all of whose eigenvalues have norm larger than $1$.
\end{definition}

\begin{definition}[Example.]
Free nilpotent groups and Heisenberg groups are expanding.
\end{definition}

The first non-expanding torsion free nilpotent group  was implicitly produced by Dixmier and Lister in \cite{dixmier1957derivations}. Later,  Dyer \cite{dyer1970nilpotent} constructed a unipotent group $U$ with a unipotent automorphism group $\Aut(U\kern 0.05 em)$. Our interest in Dyer's counterexample stems from the fact that, being unipotent, $\Aut(U\kern 0.05 em)\subset \Hom(U,U\kern 0.05 em)$ is a closed subset. Since it is also an open subset, we can conclude that $\Hom(U,U\kern 0.05 em)$ is a disconnected topological space.  Moreover, since the Lie algebra of $U$ has rational structure constants, there is a (nilpotent) lattice $\Gamma\subset U$  and the results of Mal'cev mentioned above ensure that $\Hom(\kern 0.05em\Gamma,U\kern 0.05em)=\Hom(U,U\kern 0.05 em)$. However, since the maximal compact subgroup $K$ of a unipotent group $U$ is always trivial, $\Hom(\kern 0.05em\Gamma,K\kern 0.05em)$ consists of a single point. In particular, $\Hom(\kern 0.05em\Gamma,K\kern 0.05em)$ is not homotopy equivalent to $\Hom(\kern 0.05em\Gamma,U\kern 0.05em)$.

On the other hand, when $\Gamma$ is an expanding nilpotent group the situation is entirely different. If we let $U$ be the Mal'cev completion of $\Gamma$ and $U'$ be any unipotent Lie group, the representation space $\Hom(\kern 0.05 em\Gamma,U')=\Hom(U,U')$ is contractible. Indeed, since $\Aut(U\kern 0.05 em)$ is an algebraic group, it has finitely many connected components. If $\sigma$ is an expanding automorphism of $\Gamma$, it follows that  $\sigma^l\in \Aut(U\kern 0.05 em)^o$ for some positive integer $l$. We can now choose a one-parameter-subgroup of automorphisms of $\Gamma$ containing $\sigma^l$ and reversing the associated flow induces the desired contraction.  This is the observation that allowed Silberman and Souto to prove \hyperlink{Theorem I}{Theorem I} for expanding nilpotent groups. 

%
%
%
\bibliographystyle{gtart}
\bibliography{BergeronNilpotent}

\begin{thebibliography}{}
\providecommand\bibmarginpar{\leavevmode\marginpar}
\def\urlstyle#1{{\tt #1}}

\bibitem{adem2007commuting}
\textbf{A Adem}, \textbf{F\,R Cohen}, \emph{Commuting elements and spaces of
  homomorphisms}, Mathematische Annalen 338 (2007) 587--626

\bibitem{borel1991linear}
\textbf{A Borel}, \emph{Linear algebraic groups}, volume 126, Springer-Verlag,
  New York (1991)

\bibitem{borel2002almost}
\textbf{A Borel}, \textbf{R\,S Friedman}, \textbf{J\,W Morgan}, \emph{Almost
  commuting elements in compact Lie groups}, 747, Amer Mathematical Society
  (2002)

\bibitem{brion2010introduction}
\textbf{M Brion}, \emph{Introduction to actions of algebraic groups}, Les cours
  du CIRM 1 (2010) 1--22

\bibitem{dixmier1957derivations}
\textbf{J Dixmier}, \textbf{W\,G Lister}, \emph{Derivations of nilpotent Lie
  algebras}, Proceedings of the American Mathematical Society 8 (1957) 155--158

\bibitem{dyer1970nilpotent}
\textbf{J\,L Dyer}, \emph{A nilpotent Lie algebra with nilpotent automorphism
  group}, Bulletin of the American Mathematical Society 76 (1970) 52--56

\bibitem{florentino2009topology}
\textbf{C Florentino}, \textbf{S Lawton}, \emph{The topology of moduli spaces
  of free group representations}, Mathematische Annalen 345 (2009) 453--489

\bibitem{florentino2014topology}
\textbf{C Florentino}, \textbf{S Lawton}, \emph{Topology of character varieties
  of abelian groups}, Topology and its Applications 173 (2014) 32--58

\bibitem{gomez2012fundamental}
\textbf{J\,M G{\'o}mez}, \textbf{A Pettet}, \textbf{J Souto}, \emph{On the
  fundamental group of \emph{Hom}$(\mathbb{Z}^k,{G})$}, Mathematische
  Zeitschrift 271 (2012) 33--44

\bibitem{grone1987normal}
\textbf{R Grone}, \textbf{C\,R Johnson}, \textbf{E\,M Sa}, \textbf{H
  Wolkowicz}, \emph{Normal matrices}, Linear Algebra and its Applications 87
  (1987) 213--225

\bibitem{horn1990matrix}
\textbf{R\,A Horn}, \textbf{C\,R Johnson}, \emph{Matrix analysis}, Cambridge
  university press (1990)

\bibitem{jeffrey2005flat}
\textbf{L\,C Jeffrey}, \emph{Flat connections on oriented 2-manifolds},
  Bulletin of the London Mathematical Society 37 (2005) 1--14

\bibitem{kac2000vacuum}
\textbf{V\,G Kac}, \textbf{A\,V Smilga}, et~al., \emph{Vacuum structure in
  supersymmetric Yang-Mills theories with any gauge group}, The many Faces of
  the Superworld  (2000) 185--234

\bibitem{kempf1979length}
\textbf{G Kempf}, \textbf{L Ness}, \emph{The length of vectors in
  representation spaces}, Algebraic geometry  (1979) 233--243

\bibitem{lojasiewicz1965ensembles}
\textbf{S Lojasiewicz}, \emph{Ensembles semi-analytiques}, Universit{\'e} de
  Gracovie (1965)

\bibitem{lubotzky1985varieties}
\textbf{A Lubotzky}, \textbf{A\,R Magid}, \emph{Varieties of representations of
  finitely generated groups}, 58-336, Amer Mathematical Society (1985)

\bibitem{luna1975certaines}
\textbf{D Luna}, \emph{Sur certaines op{\'e}rations diff{\'e}rentiables des
  groupes de Lie}, American Journal of Mathematics 97 (1975) 172--181

\bibitem{luna1976fonctions}
\textbf{D Luna}, \emph{Fonctions diff{\'e}rentiables invariantes sous
  l'op{\'e}ration d'un groupe r{\'e}ductif}, from: ``Annales de l'institut
  Fourier'', volume~26, Institut Fourier (1976)  33--49

\bibitem{mostow1955self}
\textbf{G\,D Mostow}, \emph{Self-adjoint groups}, The Annals of Mathematics 62
  (1955) 44--55

\bibitem{narasimhan1964holomorphic}
\textbf{M\,S Narasimhan}, \textbf{C\,S Seshadri}, \emph{Holomorphic vector
  bundles on a compact Riemann surface}, Mathematische Annalen 155 (1964)
  69--80

\bibitem{neeman1985topology}
\textbf{A Neeman}, \emph{The topology of quotient varieties}, The Annals of
  Mathematics 122 (1985) 419--459

\bibitem{onishchik1990lie}
\textbf{A\,L Onishchik}, \textbf{E\,B Vinberg}, \emph{Lie groups and algebraic
  groups}, volume 212, Springer-Verlag Berlin (1990)

\bibitem{pettet2011commuting}
\textbf{A Pettet}, \textbf{J Souto}, \emph{Commuting tuples in reductive groups
  and their maximal compact subgroups}, Geometry \& Topology 17 (2013)
  2513--2593

\bibitem{raghunathan1972discrete}
\textbf{M\,S Raghunathan}, \emph{Discrete subgroups of Lie groups}, volume~55,
  Springer-Verlag Berlin (1972)

\bibitem{richardson1988conjugacy}
\textbf{R\,W Richardson}, \emph{Conjugacy classes of n-tuples in Lie algebras
  and algebraic groups}, Duke Math. J 57 (1988) 1--35

\bibitem{richardson1990minimum}
\textbf{R\,W Richardson}, \textbf{P\,J Slodowy}, \emph{Minimum vectors for real
  reductive algebraic groups}, Journal of the London Mathematical Society 2
  (1990) 409--429

\bibitem{schwarz1989topology}
\textbf{G\,W Schwarz}, \emph{The topology of algebraic quotients}, Topological
  methods in algebraic transformation groups  (1989) 135--151

\bibitem{sikora2012character}
\textbf{A Sikora}, \emph{Character varieties}, Transactions of the American
  Mathematical Society 364 (2012) 5173--5208

\bibitem{simpson1994moduli}
\textbf{C\,T Simpson}, \emph{Moduli of representations of the fundamental group
  of a smooth projective variety I}, Publications math{\'e}matiques de l'IHES
  79 (1994) 47--129

\bibitem{witten1982constraints}
\textbf{E Witten}, \emph{Constraints on supersymmetry breaking}, Nuclear
  Physics B 202 (1982) 253--316

\bibitem{witten1998toroidal}
\textbf{E Witten}, \emph{Toroidal compactification without vector structure},
  Journal of High Energy Physics 1998 (1998) 006

\end{thebibliography}
\end{document}